\documentclass[12pt,leqno]{article}
\usepackage{amsmath,amsfonts,amscd,amssymb,latexsym}
\usepackage{amsfonts,amsmath,latexsym}
\newtheorem{lema}{Lemma}
\newtheorem{teo}[lema]{Theorem}

\newtheorem{propos}[lema]{Proposition}

\newtheorem{coro}[lema]{Corollary}
\newtheorem{defi}[lema]{Definition}
\newtheorem{rema}[lema]{Remark}

\newcommand{\br}{\mathbb R}

\newcommand{\bz}{\mathbb Z}

\newcommand{\spe}{{\rm {sp\,}}}

\newcommand{\proof}{\noindent{\sc {Proof:\,}}}

\newcommand{\hbx}{\hfill$\Box$}
\def \endproof{\hbx}

\begin{document}

\title{Estimates of the topological entropy from below for continuous self-maps on
some compact manifolds}

{\author{Wac{\l}aw Marzantowicz \& Feliks Przytycki }}
\date{\today}

\maketitle

\begin{abstract}
 \noindent
Extending our results of \cite{MrPr}, we confirm that Entropy
Conjecture holds for every continuous self-map of a compact
$K(\pi,1)$ manifold with the fundamental group $\pi$ torsion free
and virtually nilpotent, in particular for  every continuous map
of an infra-nilmanifold. In fact we prove a stronger version,
a lower estimate of the topological entropy of a map by logarithm of 
the spectral radius of an associated  "linearization matrix" with integer entries.

From this, referring to known estimates of Mahler measure of
polynomials, we deduce some absolute lower bounds for the entropy.

\end{abstract}

\section{Introduction}\label{covers}

Let $M$ be a compact manifold and $f:M \to M$ a continuous
self-map of $M$. The topological entropy $h_{\rm{top}}(f)$,
denoted shortly by ${\bf h}(f)$, is defined as $\lim_{\epsilon\to
0}\limsup_{n\to\infty} 1/n \log \sup \#Q$, supremum over all $Q$
being $(\epsilon, n)$-separated. $Q$ is called $(\epsilon,
n)$-separated, if for every two distinct points $x,y\in Q$,
$\max_{j=0,...,n}d(f^j(x),f^j(y))\ge \epsilon$. Here
 $d$ is a metric on $M$ consistent with the topology; in fact
${\bf h}(f)$ does not depend on the metric (cf. \cite{Szl}).

Entropy Conjecture,  denoted shortly as EC,  says that the
topological entropy of $f$ is larger or equal to the logarithm of
the spectral radius of the linear operator induced by $f$ on the linear
spaces of  cohomology of $M$ with real coefficients. It was posed
by M. Shub in seventies who asked what suppositions on $f$ and $M$
imply EC.

\

We prove the following

\

\noindent{\bf Theorem A} {\it Assume that  a compact  manifold $M$
is a $K(\pi,1)$-space  with the fundamental group $\pi$ being
torsion free and virtually nilpotent. Then  {\rm EC} holds  for
every continuous self-map $f$ of $M$. In particular {\rm
EC} holds for every continuous self-map of any compact
infra-nilmanifold.}

\

$\bullet$ A group $\pi$ is called {\it virtually nilpotent} if it contains a
finite index nilpotent subgroup $\Gamma$.  We can assume that
$\Gamma$ is a normal subgroup of $\pi$.

(Indeed, for any pair of groups $K\subset L$, where $K$ has a
finite index in $L$, one has a homomorphism $\rho: L \to {\rm
Sym}(L/K)$ into the symmetry group of the quotient space. Then
$\ker \rho$ is a normal subgroup of $L$, it has finite index, and
 is contained in  $K$.)

$\bullet$ One can replace the assumption that $\pi$ is virtually nilpotent by the assumption
that $\pi$ has polynomial growth, \cite{Gr1}.

$\bullet$ Theorem A is a step towards proving a conjecture by A. Katok \cite{Ka}
saying that EC holds for every
continuous map if the universal cover of $M$ is homeomorphic to an
Euclidean space $\br^d$.


$\bullet$ One can even ask whether EC holds for
every continuous self-map of a $K(\pi,1))$ compact  manifold (or a finite 
CW-complex).

\


{\bf 1.1. Affine maps}

\

We refer to the following theorem by A. Malcev and L. Auslander
about the {\it existence of a model} \cite{Franks} p.76:

Assume that $\pi$ is finitely generated torsion free virtually nilpotent group.
Then it contains a finite index normal subgroup $\Gamma$ which can
be embedded as a lattice, i.e. a discrete co-compact subgroup,  in
a connected, simply connected nilpotent Lie group $G$. The embedding can be
extended to the embedding of
$\pi$  in the group $\rm{Aff}(G)$ of affine
mappings of $G$, so that $\pi \cap G = \Gamma$. More precisely, if $
\mathbf{C}\subset {\rm Aut}(G)$ denotes the maximal compact
subgroup of the group of automorphisms of $G$, then $\pi \subset
G \ltimes{\mathbf{C}}\subset G \ltimes{\rm Aut}(G)$; this
embedding of $\pi$ is called an {\it almost Bieberbach group}.

It follows then from the definition that $\pi$ acts on $G$
properly discontinuously (First note that
if $\alpha\in \pi$ has a fixed point
$z\in G$, then $\alpha^\ell(z)=z$ and $\alpha^\ell\in \Gamma$ for $\ell= \#(\pi/\Gamma)$.
Then $\alpha^\ell=e$ the unity of $\pi$, hence $\alpha=e$ by the assumption that $\pi$ is torsion free.)

The quotient manifold $\mathbf{IN}=G/\pi$ is called an {\it
infra-nilmanifold}. It is regularly finitely covered by the
nilmanifold $\mathbf{N}=G/\Gamma$, with the deck transformation
group equal to $H=\pi/\Gamma$.

Note that every compact manifold finitely covered by a nilmanifold, in particular
every infra-nilmanifold, satisfies the assumptions of Theorem A.
Indeed,
$G$ is homeomorphic to $\br^d$, where $d=\dim M$ (cf. \cite{Ra}).
If $\pi$ were not torsion free, a cyclic subgroup
generated by an element $\{g\} \simeq \bz_p$ of a prime order
would act freely on $\br^d$. The latter is impossible, as follows
from the Smith theory (cf. \cite{Br}).

The image of  the embedding of  $\pi$ into $ G \ltimes \mathbf{C}\subset {\rm Aut}(G)$ will
be denoted by $\pi_{\mathbf{IN}}$. It is a deck transformation group of the cover  
$p_{\mathbf{IN}} :G\to G/\pi_{\mathbf{IN}}$ and, distinguishing an arbitrary $z\in G$, 
it can be identified in a
standard way with the fundamental group  $\pi_1({\mathbf{IN}},
p_{\mathbf{IN}}(z))$
of ${\mathbf{IN}}$.

Consider now any $M$ being $K(\pi,1)$ as in Theorem A. 
The group $\pi$ acts properly discontinuously on the universal cover space $\tilde M$ 
and it can be identified with $\pi_M$, which is the deck transformation group  
of the universal cover $p_M:\tilde M\to M$.  
 
(Similarly to ${\mathbf{N}}\to {\mathbf{IN}}$, we have a regular 
finite cover  $\tilde M/\Gamma\to M=\tilde M/\pi$, with the deck transformation group $H$.)

Every continuous $f:M\to M$ induces an endomorphism $F=F_f$ of $\pi_M$, unique up to 
an inner authomorphism.  To define this use 
$f_\#:\pi_1(M,z)\to\pi_1(M,f(z))$
between the fundamental groups and standard identifications of these groups with $\pi_M$. 
For more details see Section 3.

When we identify $\pi$ with $\pi_{\mathbf{IN}}$, then we can consider $F$ 
as an endomorphism of $\pi=\pi_{\mathbf{IN}}$. 

By K. B. Lee, \cite[Th.1.1 and Cor 1.2]{LeeK}, there exists an
affine self-map $\Phi=\Phi_f=(b,B)$ of $G$, with $b\in G$, $B\in
{\rm End}(G)$, such that  for all $x\in G, \alpha\in
\pi_{\mathbf{IN}}$ we have
$$ \Phi(\alpha(x))   = F(\alpha)( \Phi(x)).  \eqno (1)
$$

Hence, in view of (1), there exists a factor $\phi=\phi_f$ of $\Phi$ on
$\mathbf{IN}$ under the action of $\pi_{\mathbf{IN}}$. In general
one calls factors of affine $\Phi$ on $G$ satisfying (1), affine
maps on $\mathbf{IN}$, in particular  $\phi_f$ is an affine map on
$\mathbf{IN}$.

Since both $M$ and $\mathbf{IN}$ are $K(\pi, 1)$-spaces, they are
homotopically equivalent, \cite{Sp}. 
One can find a homotopy equivalence 
$h:M\to{\mathbf{IN}}$ inducing our identification between $\pi_M$ and $\pi_{\mathbf{IN}}$.
Then $\phi\circ h$ is homotopic to $h\circ f$, see \cite[Section 8.1]{Sp} and Section 3.   
Note, \cite{LeeK}, that if $M=\mathbf{IN}$, i.e. $f$ is a
self-map of an infra-nilmanifold $\mathbf{IN}$, then $\phi_f$ is
homotopic to $f$.

Let $\Phi=\Phi_f=(b,B)$ be an affine self-map of $G$ associate
to $f:M\to M$ (in fact to the homomorphism $F$).
The differential at the unity, $D(B)(e)$, is an
endomorphism of $\cal G$, the Lie algebra of $G$. Denote $D(B)(e)$
by $D_f$. Let $\sigma(D_f)=\{\lambda_1,\dots ,\lambda_d\}$ be the
set of all its eigenvalues counted with multiplicities. Consider
$\spe(\wedge^* D_f)$, the spectral radius of the full exterior
power $\wedge^*D_f={\underset{0}{\overset{d}\oplus}} \wedge^k D_f
$ of $D_f$. Of course it is equal to $\prod_{j:|\lambda_j|>1} |\lambda_j|$, where
$\lambda_j\in \sigma(D_f)$, provided at least one $\lambda_j$ has absolute value larger than 1.
Otherwise it is equal to 1 (in $\wedge^0 D_f$).

\

{\bf 1.2. Linearization matrices}

\



One can assign to an endomorphism $F=f_\#$ not only a linear map
$D_f: \br^d \to \br^d$, ${\cal G}\equiv {{\br}}^d$, but also an
integer $d\times d$ matrix $A_{[f]}$ called the linearization of
the homotopy class $[f]$, because  endomorphisms $f_\#$ are in the
one-one correspondence with homotopy classes $[f]$ of self-maps of
a $K(\pi,1)$-space. Sometimes we shall use the notation $A_f$.

An endomorphism $F: \pi \to \pi$ does not preserve the nilpotent
subgroup $\Gamma \lhd \pi$ in general. But $\Gamma$ contains a
subgroup $\Gamma^\prime \lhd \pi$ such that $\Gamma^\prime$ is
nilpotent has finite index in $\pi$ and is invariant for $F$ (cf.
\cite{LeeJLeeK} or Proposition \ref{existence of admissible}).
Since $\Gamma^\prime$ has finite index in $\Gamma$ it is also a
lattice in $G$.
By (1) the endomorphism $B:G \to G$ is an
extension of $F: \Gamma^\prime \to \Gamma^\prime$, see Section 2.

 Let $\,G\,$ be a nilpotent connected simple-connected  Lie group, $\,\Gamma\,$ its
  lattice.   By the definition the descending central series of commutators
$\,G_0=G\,,$ $\,G_{i+1}=[G,G_i]\,$ of $G$ is finite  with the last
group trivial $\,G_k=\{e\}\,.$
 For each $i$ define $\Gamma_i=G_i\cap \Gamma$.   Each $\Gamma_i$ is a lattice in $G_i$,
as follows from Malcev theorem \cite[Theorem 1]{Malcev}.
Consequently, if
    $B: G\to G$ is a homomorphism of $G$ preserving a lattice $\Gamma$ it preserves
    each subgroup $\Gamma_i$ thus it induces an endomorphism $B_i$ on each factor group
    $\Gamma_i/\Gamma_{i+1}$, $0\leq i \leq k$. Clearly
    $\Gamma_i/\Gamma_{i+1}$ is abelian and
    torsion free of dimension $d_i$.
Therefore the action of $B_i$ on
$\Gamma_i/\Gamma_{i+1}\equiv
    \bz^{d_i}$ yields an integer  $d_i\times d_i$ matrix  $A_i$
    which is uniquely defined up to a choice of basis, i.e. up to a
    conjugation by a unimodular matrix.
Finally we put $A_{[f]}:=
    {\underset{i=0}{\overset{k}\oplus}} A_i$, which defines an
   integer $d\times d$ matrix, because
  ${\underset{i=0}{\overset{k}\oplus}} d_i = d$.  (See \cite{KMC} and \cite{JeMr2} for some
extensions to solvmanifolds.)

It is easy to show that $\sigma(A_{[f]})= \sigma(DB(e))$, thus
consequently $\sigma(\wedge^*A_{[f]})= \sigma(\wedge^*D_f)$, which
implies $\spe(\wedge^*A_{[f]})= \spe(\wedge^*D_f)$ (cf.
\cite{MrPr}, \cite{JeMr2}).



In fact the matrices $A$ can be defined directly,
using $F=f_\#:\pi\to\pi$,
without constructing $G$. One can use the  series of {\it isolators}
$\sqrt[\Gamma]{\Gamma'_i}=\{x\in \Gamma: (\exists\ell>0) x^\ell\in\Gamma'_i\}$,
for $\Gamma'_i$ being commutators in the descending central series for $\Gamma$,
 i.e.  $\Gamma'_{i+1}=[\Gamma,\Gamma'_i]$, see \cite{Dekim1} Section 3 and
for example \cite[Lemma 11.1.8]{Passman}. In fact
$\sqrt[\Gamma]{\Gamma'_i}=\Gamma_i=G_i\cap \Gamma$ defined above, see 
\cite[Lemma 1.2.6]{Dekim}.

\

{\bf 1.3. Conclusion}

We are in the position to formulate a sharper version of Theorem A,
namely

\

\noindent {\bf Theorem B} {\it For every continuous self-map $f$
of a compact  manifold $M$  which is a  $K(\pi,1)$-space  with the
group $\pi$ torsion free and virtually nilpotent,

$$
{\bf h}(f)\,\geq \, \log \spe(\wedge^*
 D_f)\, = \, \log \spe(\wedge^* A_{[f]})\,.
$$
 In the case $M$ is an infra-nilmanifold the
equality holds for
every
affine map
$\phi: M\to M$, a factor of an affine $\Phi$ satisfying (1), in particular for $\phi_f$.
In consequence
for every continuous self-map $f$ of $M$ we have  ${\bf h}(f)\geq {\bf
h}(\phi_f)$, i.e. $\phi_f$ minimizes entropy in the homotopy class
of $f$.}

\vskip0.13cm

Maybe considering of $M$ is not needed and it is sufficient to
consider only $\mathbf{IN}$. Since  the topological entropy is an
invariant of conjugation  by a homeomorphism,  this would follow
from Borel conjecture, which states that the fundamental group
$\pi$ of a manifold being $K(\pi,1)$-space determines $M$ up to
homeomorphism,  This has been confirmed by Farrell and Jones
\cite{Fa-Jo} for a class of groups that contains the
almost-Bieberbach  groups, except in dimension 3.

Formulating Theorems A and B we have followed a suggestion  by M. Shub
\cite{Sh2}, to assume a discrete group point of view. Having given
an endomorphism $F=f_\#:\pi\to\pi$ of a finitely generated torsion free virtually
nilpotent group, we associate to it a linear operator $D_f$, or an
integer $d\times d$  matrix $A_f$. As suggested in \cite{Sh2} the
logarithm of  spectral radius of $\spe(\wedge^*D_f)$, or
$\spe(\wedge^*A_f)$,  is "a kind of volume growth" of $f_\#$. In
Theorem A it is replaced by the spectral radius of the map induced
on real cohomologies of the group $\pi$.

\

We shall present two proofs of Theorems A and B.

The first one, in Section 2,  concludes Theorems A and B from
analogous theorems in \cite{MrPr}, with $f:\mathbf{N}\to
\mathbf{N}$ a continuous  map of a compact nilmanifold. However
this proof of Theorems A and B does not work in dimension 3.

The second proof, in Section 3, holds for $f$ on $M$ and uses only
a homotopy equivalence between $M$ and $\mathbf{IN}$. It directly
repeats the arguments of \cite{MrPr}.

\

An important observation is that $A_{[f]}$ is an integer matrix. This allows, in Section 4,
to prove absolute estimates from below  for
$\spe(\wedge^* A_f)$, where $f$ is an expanding map of a compact
manifold (without boundary) or  an Anosov diffeomorphism of a
compact infra-nilmanifold. The latter uses  number theory results
estimating the Mahler measure of an integer polynomial.

\

The authors would like to express their thanks  to K. Dekimpe, E.
Dobrowolski,  T. Farrell, A. Katok and M. Shub for helpful conversation.

\

\section{Entropy Conjecture on infra-nilmanifolds}\label{EC}

 The proof of Theorems A and B we present in this section
will hold for $M$ finitely
covered by a nilmanifold and follows from two standard facts and the
main theorem of \cite{MrPr} in which the topological entropy of a
continuous map of nilmanifold is estimated by the corresponding
quantities. We begin with the following

\begin{propos}\label{entropy for covering} Suppose that we have finite
cover $(\tilde{M}, p, M)$ of compact metric spaces, i.e.
$\tilde{M}$ is the total space of the cover, $M$ the base space,
and $p:\tilde{M}\to M$ the covering space. Let a pair
$(\tilde{f},f)$, $\tilde{f}:\tilde{M} \to \tilde{M}$, $f:M\to M$,
be a map of a this covering, i.e. $p \tilde{f} = f p$. Then $$
{\bf h}(f) ={\bf h}(\tilde{f})\,.$$
\end{propos}

\proof It is elementary and  is given in \cite{FaSh}. Briefly: The
$p$-preimage of an $(n,\epsilon)-f$-separated set in $M$ is
$(n,\epsilon)-\tilde f$-separated (in a metric $\tilde d$ on
$\tilde M$ is being the lift of a metric $d$ on $M$ chosen to
define the entropy), hence ${\bf h}(f) \le {\bf h}(\tilde{f})\,.$
(In fact only the continuity of $p$ was substantial in this
proof).

Conversely, let $Q$ be an $(n,\epsilon)-\tilde f$-separated set in
$\tilde M$ consisting of points in a ball $B(z,\epsilon/2)$. Let
$\delta>0$ be a constant such that $p$ is injective on every ball
in $\tilde M$ of radius $\delta$. We prove that the set $p(Q)$ is
$(n,\epsilon)-f$-separated.  Indeed, take $\epsilon < \delta$ and
suppose that for $x,y\in Q$ we have
$d(f^j(p(x)),f^j(p(y))<\epsilon$ for all $j=0,1,...,n$. Let
$j_0\ge 0$ be the smallest $j\le n$ such that
$d(\tilde{f}^j(x),\tilde{f}^j(y)\ge \epsilon$. Then
$d(\tilde{f}^j(x),\tilde{f}^j(y)\ge \delta-\epsilon$ (i.e.
projections by $p$ are close to each other but the points are in
different components of preimages of a small ball under the cover
map). This is not possible for $j=0$ by $Q\subset
B(z,\epsilon/2)$. If it happens for another $j$ it means that the
$\tilde f$ image of two points within  the distance $< \epsilon$
have distance $\ge \delta-\epsilon$, what for $\epsilon$ small
enough contradicts the uniform continuity of $\tilde f$.
\endproof

\

The same trick as above has been used in \cite[Proof of
Theorem 2.2]{MrPr}, to which we refer here in Section 3,  in the setting
of the projection $p$ from the universal cover of $M$ to $M$.

\begin{propos}\label{cohomology of
regular cover} Let $(\tilde{M}, p, M)$ be a regular covering over a finite CW-complex, e.g. over a compact
manifold with the deck transformation group $H$. Let next  $f:M\to M $ be a continuous map   and
$\tilde{f}:\tilde{M} \to \tilde{M}$ its lift to $\tilde{M}$.

Then $\sigma(H^*(f))\subset \sigma(H^*(\tilde{f})) $ and consequently ${\rm sp} (H^*(f)) \leq {\rm sp}
(H^*(\tilde{f}))$.
\end{propos}

\proof  The action of $H$ on $\tilde{M}$ defines an action on the
real cohomology space for any $g\in H$ given by $g \mapsto \,
H^*(g) :H^*(\tilde{M})\to H^*(\tilde{M})$.
 Since $M$ is the orbit space of a free action of a
finite group $H$ on $\tilde{M}$, we have  $H^*(M;\br) \simeq
H^*(\tilde{M};\br)^H$, when the latter is the fixed point of
action of $H$ on the real cohomology spaces, i.e. the image of  a
linear projection $\frac{1}{\vert G \vert}\sum_{g\in H} H^*(g)$
(cf. \cite{Br} Section III.2). Moreover this isomorphism is given by the map
$H^*(p) :H^*(M;\br) \to H^*(\tilde{M};\br)$ induced by the
covering  projection $p$. Since $p\,\tilde{f} = f \,p $, we have
$H^*(\tilde{f})H^*(p) =H^*(p) H^*(f)$. This means that the linear
subspace $H^*(M;\br)= H^*(\tilde{M};\br)^H ={\rm im}\, H^*(p)$ is
preserved by $H^*(\tilde{f})$. Consequently $H^*(f)$ can be
identified with a restriction of a linear map $H^*(\tilde{f})$ to
an invariant linear subspace, which shows that $\sigma(H^*(f))
\subset \sigma(H^*(\tilde{f}))$, and consequently proves the
statement. (Note that $\tilde{f}$ is not an $H$-equivariant map in
general.)
\endproof

\

An immediate corollary is

\begin{coro}\label{ECcover} Let $(\tilde{M}, p, M)$ be a regular cover over a finite CW-complex,
e.g. over a compact
manifold, $f:M\to M$ is a continuous map which has a lift $\tilde f$ to $\tilde M$ and EC holds for
$\tilde f$, then EC holds for $f$.
\end{coro}

A problem in a general situation is with the existence of the lift $\tilde f$. Fortunately one can go around it for $\tilde M$ being a nilmanifold.

\begin{defi}\label{admissible group}
Let $\Gamma\triangleleft\pi$ be a normal nilpotent subgroup of finite index in
$\pi$ and let let $s$ be an endomorphism of $\pi$. We
say that a group $\Gamma^\prime \subset \Gamma$ is $s$
-admissible if
\begin{itemize}
\item[1)] {$s(\Gamma^\prime)\subset \Gamma^\prime,$} \item[2)]
{$\Gamma^\prime \subset\Gamma$ is normal in $\pi$ and
$[\Gamma:\Gamma^\prime] <\infty$, i.e. $\Gamma^\prime $ has a
finite index in $\Gamma$.}
\end{itemize}
\end{defi}



\begin{propos}\label{existence of
admissible}

For a nilpotent group $\Gamma$ normal and of finite index in a group $\pi$
there exists
a group $\Gamma'\subset\Gamma$, \ $\Gamma'\triangleleft\pi$,
admissible for every endomorphism $s$ of $\pi$.
(sometimes such $\Gamma'$ is called a fully characteristic subgroup).
\end{propos}

\proof Repeat verbatim the argument of Lemma 3.1 of \cite{LeeJLeeK} and define
$$     \Gamma^\prime= {\text{group generated by}}\;\;
\{ \gamma^k: \, \gamma \in \pi\},\;\; {\text{where}}\;\, k\;\,{\text{is the order of}}\;\, H=\Gamma/\Gamma'\,.$$
 This is a subgroup preserved by every endomorphism of $\pi$, in particular
$\pi$ is normal.

Next we define a group
 $$\Gamma(k):= \;{\text{group generated by}}\;\; \{x^k\}:
x\in \Gamma\,.$$ Of course $\Gamma(k)\subset \Gamma'$. It is
enough to show that $\Gamma(k)$ is of a finite index in $\Gamma$.
Apply an argument used in \cite{LeeJLeeK}: Since $\Gamma$ is
nilpotent, it is polycyclic, \cite{Ra}. But for any polycyclic
group $\Gamma$ the corresponding group $\Gamma(k)$ has a finite
index, cf. \cite[Lemma 4.1]{Ra}. In particular adapting the
argument of   Lemma 4.1 of \cite{Ra} one shows its
assertion for a nilpotent group, by an induction over the
 length of nilpotency.

\endproof

 Note that the number  $k$ used to define the group $\Gamma(k)$  is not unique, e.g. we can take any its
 multiple getting a smaller group with the required property.
To get a larger group $\Gamma^\prime$ than that of  \cite{LeeJLeeK}
   we can use $k$ equal to the ${\rm LCM}\{\#
h; h\in H  \}$  the order of element, instead of $k=\# H$ the order of $H$.

\begin{coro}\label{lift to nil} For any  compact manifold $M$
finitely
covered by a nilmanifold
$\mathbf{N}$
there exists a regular finite
cover $(\tilde{\mathbf{N}}, \tilde{p}, M)$ of $M$ by a nilmanifold
$\tilde{\mathbf{N}}$ such that
every continuous map
$f: M\to M$ has a lift $\tilde{f}: \tilde{\mathbf{N}}
\to \tilde{\mathbf{N}}$, i.e. $ p\, \tilde{f} = f \,p$.

\end{coro}

\proof The assertion follows from  Proposition \ref{existence of admissible}
for $\pi$ the fundamental group of $M$, \
$\mathbf{N}=G/\Gamma$ for $\Gamma$ a subgroup of $\pi$ and $s=f_\#$. We can assume that $\Gamma$ (hence $\Gamma'$) is normal in $\pi$, see Introduction.
We define $\tilde{\mathbf{N}}:=G/\Gamma'$.
A lift $\tilde f$ exists, since the homomorphism $f_\#: \pi \to \pi $ preserves  $\Gamma'$
identified with $\tilde{p}_\#\pi(\tilde{\mathbf{N}})$, see \cite{Sp}.

\endproof

\

Together with Corollary \ref{ECcover} and EC for nilmanifolds,
\cite{MrPr}, this proves EC for all continuous self maps of $M$
finitely covered by nilmanifolds, in particular for all $M$ being
infranilmanifolds.

\

\proof {\sc of Theorem B } (for $M$ finitely covered by a nilmanifold $G/\Gamma$).

As above we find an admissible $\Gamma'\subset\Gamma$.
The quality (1)
for $\alpha=g\in\Gamma'$  takes the form
$$
B(g(x))b=F(g)B(x)b.
$$
As the left hand side expression is equal to $B(g)B(x)b$, we get
$F(g)=B(g)$, i.e. $f_\#|_{\Gamma'}=B|_{\Gamma'}$.

(In fact $B$ in $(b,B)$ was found in
\cite{LeeK} just as an extension of $f_\#|\Gamma'$ to $G$.)

$B$ is a lift to $G/\Gamma'$ of $\phi_f$ homotopic to $f$.
Theorem B follows from Corollary \ref{lift to nil} ($B$ is $\tilde f$ there), from Proposition
\ref{entropy for covering}, and from Theorem B for self-maps $f$ of nilmanifolds, \cite{MrPr}.

\endproof





\







\vskip 0.0cm

\section{Another proof of EC}\label{EC2}

Now we provide another proof of Theorems A and B without additional assumptions
by showing that a modification of the proof for the
nilmanifolds  given in \cite{MrPr} works.

\

Let us remind the notation. 
$M$ is a compact manifold, being $K(\pi,1)$ for $\pi$ a virtually nilpotent 
torsion free group $\pi$.\ $G$ is a connected
simply connected nilpotent Lie group and ${\mathbf{IN}}=G/\pi$ where $\pi$ is 
embedded in $\rm{Aff}(G)$ as $\pi_{\mathbf{IN}}$,
acting discontinuously on $G$ so that ${\mathbf{IN}}$ is an infra-nilmanifold,
see Existence of a Model Theorem in  Introduction.
We have the universal covers $p_M:\tilde M\to M$ and $p_{\mathbf{IN}}:G\to {\mathbf{IN}}$.

Remark that we use the  right action, thus ${\mathbf{IN}}=G/\pi$, instead for 
$\pi\backslash G$ used in \cite{LeeK}. Then the action of an
affine map $(d, D)$, $d\in G$, $D\in {\rm Endo}(G)$, is given as
$(d,D)x= (Dx)d$.

\ We assume that all  metrics under consideration are induced by
Riemannian metrics. We need the following

\begin{lema}\label{equivalence of metrics} {\rm [Lemma on the equivalence of metrics on $G$]} Any right
invariant metric $\rho$ on $G$ is equivalent to $\tau_G$ (i.e. the
mutual ratios are bounded)  being a lift of an arbitrary metric
$\tau_{\mathbf{IN}}$ on an infra-nilmanifold ${\mathbf{IN}}=G/\pi$.

\end{lema}

\proof By compactness the lift $\tau_\Gamma$ of $\tau_{\mathbf{IN}}$ to
$G/\Gamma$, where $\Gamma= G\cap \pi$, is equivalent to
$\rho_\Gamma$, the projection of $\rho$ to $G/\Gamma$. Therefore
the lifts to $G$ are also equivalent.\endproof


%

Let us stop for a while on the homotopy equivalence between $M$ and ${\mathbf{IN}}$ 
making some explanations from Introduction more precise.

\begin{lema}\label{homotopy}  
There exists a homotopy equivalence $h:M\to \mathbf{IN}$. Moreover for 
every continuous $f:M\to M$ there exists an affine mapping $\phi: \mathbf{IN}\to \mathbf{IN}$ 
such that $\phi\circ h \simeq h\circ f$, where $\simeq$ means: homotopic. 
\end{lema}

\proof  The action $\pi_M$ of $\pi$ on $\tilde M$ can be identified, chosen 
an arbitrary distinguished 
point $\tilde z_0\in \tilde M$,  with the fundamental group $\pi_1(M,p_M(\tilde z_0))$ 
by projecting by $p_M$ of curves joining $\pi_M(\tilde z_0)$ to $\tilde z_0$.
Similarly $\pi_{\mathbf{IN}}$ can be identified with 
$\pi_1({\mathbf{IN}}, p_{\mathbf{IN}}(\tilde w_0)$ 
for an arbitrary distinguished point $\tilde w_0\in G$. 

Let $h:M\to  {\mathbf{IN}}$ be a homotopy equivalence; 
its existence follows from \cite[Section 8.1]{Sp}. 
Choose its lift $\tilde h:\tilde M \to G$ and distinguish $\tilde z_0,\tilde w_0$ such that 
$\tilde h(\tilde z_0)=\tilde w_0$. Denote $z_0=p_M(\tilde z_0)$ and $w_0=p_{\mathbf{IN}}(\tilde w_0)$. 
We have $h_\#:\pi_1(M,z_0)\to\pi_1({\mathbf{IN}},w_0)$ which yields, with respect to distinguished 
$\tilde z_0, \tilde w_0$, an isomorphism between the deck transformation groups 
$H:\pi_M\to\pi_{\mathbf{IN}}$.

Let $\tilde f: \tilde M\to\tilde M$ be a lift of $f$. 
Define $F_M: \pi_M\to\pi_M$, with respect to $\tilde z_0$ and $\tilde f (z_0)$ (similarly to the way we defined $H$).
Finally define $F=F_{\mathbf{IN}}:=H\circ F_M \circ H^{-1}$ and affine $\Phi$ and $\phi$ as 
in Introduction, relying on \cite{LeeK}.  

By construction we get $F_{\mathbf{IN}}\circ H =  H \circ F_M$. 
Denote $x_1=\Phi (\tilde h (\tilde z_0))$ and  
$x_2=\tilde h ( \tilde f (\tilde z_0))$. If $x_1=x_2$ then $\phi_\#h_\#=h_\#f_\#$ on $\pi_1(M,z_0)$, hence 
 $\phi\circ h \simeq h\circ f$, by \cite[Section 8.1, Theorem 11]{Sp}. 
Otherwise one can consider $k:{\mathbf{IN}}\to  {\mathbf{IN}}$ homotopic to identity 
and its lift $\tilde k$, such that $\tilde k (x_1)= x_2$. It gives (with respect to 
distinguished $x_1, x_2$) identity on $\pi_{\mathbf{IN}}$. 
Hence
$k_\#\phi_\#h_\#=h_\#f_\#$ on $\pi_1(M,z_0)$,
hence again we deduce $\phi\circ h \simeq h\circ f$.
\endproof

\

\noindent{\sc Proof (of Theorems A and B):}  Consider metrics
$\tau_{\tilde M}, \tau_G$ on $\tilde M, G$ respectively, being
lifts to the universal covers $p_M:\tilde M\to M$ and 
$p_{\mathbf{IN}}: G\to {\mathbf{IN}}$ of
arbitrary metrics $\tau_M, \tau_{\mathbf{IN}}$ on $M,{\mathbf{IN}}$.

Let $f:M\to M$ be a continuous map and $\tilde f$ its lift to $\tilde M$. 

Let $h:M\to {\mathbf{IN}}$ be a homotopy equivalence such that $h\circ f\simeq\phi_f\circ h$.
for $h$ and $\phi=\phi_f$ as in Lemma \ref{homotopy}.
Let $\tilde h:
\tilde M \to G$ be a lift of $h$. Then the distance in $\tau_G$ between 
$\tilde h \circ \tilde f$ and $\Phi\circ \tilde h$ is bounded, by a constant $\xi_1>0$ 
(since the lifts are joined by a lift of a homotopy, up to a deck transformation, that is up to an isometry). \newline


Let $x_n,\; n=0,1,2,...$ be an $\tilde f$ trajectory.





Hence $w_n=\tilde h(x_n)$ is a
$\xi_1-\Phi$-trajectory in the metric $\tau_G$, hence, by Lemma \ref{equivalence of metrics}, a 
$\xi_2-\Phi$-trajectory
in
$\rho$, the right invariant metric on $G$.

Finally for $\Phi=(b,B)$ the sequence $w_n$ is a $\xi_3-B'$-trajectory  for
$B'
=b^{-1}Bb$, i.e.
$B'(x)=b^{-1}\cdot B(x)\cdot b$.

Indeed, $\rho(\Phi(w),B'(w)) =
\rho(B(w)b,B'(w))=\rho(bB'(w),B'(w))=\rho(b,e)$, the latter equality
by the right invariance of $\rho$. Hence $w_n$ is a $\xi_2 +\rho(b,e)$
trajectory for $B'$.

Note that the spectra of the derivatives (linearizations) of $DB(e)$ and
$DB'(e)$ coincide as these operators are conjugate.

\

Now we define a mapping $\Theta$ from $(w_n)$ to a $B'$-trajectory
in $G^u$ the unstable subgroup for $B'$ by proceeding  as in
\cite{MrPr}: First we define $w_n\mapsto \pi^u(w_n)$, the
"projection" to $G^u$, i.e. we write $w_n=g^{cs}\cdot g^u$ where
$g^{cs}\in G^{cs}$ the central stable subgroup and
$\pi^u(w_n):=g^u\in G^u$. Next $\Theta(w_n)$ is defined as the
unique $B'$-trajectory in $G^u$ subexponentially "shadowing"
$\pi^u(w_n)$. Finally, we define $\theta(x):=\Theta\circ \tilde h (x)$.

For an arbitrary $\epsilon>0$, 
for $((1+\epsilon)^j, n)-B'$-separated points in $G^u$,  $j=0,...,n$, 
(contained in a small disc), i.e. such that for some $j$ their $j$-th images under $B'$ are within the distance at least $(1+\epsilon)^j$, we choose
points $w$ in their $\Theta$ preimages (also in a small disc) and next points $x$ in
$\tilde h$-preimages, in a small disc. This is a crucial point which uses the fact that
$\tilde h$ is onto, since $|{\rm{deg}}\, h|=1$, compare Remark 4.8
in \cite{MrPr}.

If two points $p_M(x),p_M(y)$ are $(\epsilon,n)-f$-close (i.e. not separated), then so are
$x,y$. Hence $\tilde h(x)$ and $\tilde h(y)$ are 
$(\xi_4, n)$-close (with respect to $\Phi$, hence $B'$) in $\rho$ for
a constant $\xi_4$. Hence their $\Theta$ images are 
$((1+\epsilon)^j,n)-B'$-close, a contradiction.
\endproof

\

Note that we did not use an admissible group constructed in
Proposition \ref{existence of admissible}

\begin{rema}\label{Szczepanski question}\rm
 The statement of Theorem A, in a weaker form
 for  the flat manifolds,
was posed as a question  by Szczepa{\'n}ski in his article \cite{Szcz}. Earlier, a very special case of entropy
conjecture an for affine map of a compact affine manifold was proved by D. Fried and M. Shub in \cite{FrSh}.
\end{rema}

\section{Absolute estimates of entropy}



Famous Lehmer  conjecture in number theory  states that
there exists a constant $C$, called Lehmer constant, such that for
every
integer
 polynomial $w(x)=a_0x^d +a_1 x^{d-1}+\, +\cdots\,+a_d$, not being a
product of cyclotomic polynomials (all zeros being roots of 1) and $x^k$,
for the Mahler
$M(w)$ measure of $w$ have
$$
M(w):= \vert a_0\vert\, {\underset{\lambda_i}\prod} \max(1,\vert \lambda_i\vert) \geq
C\,,$$
where the product is taken over all zeros of $w(x)$.

\


There are estimates of the Mahler measure which depend on the
degree of an irreducible polynomial (the degree of an algebraic
number).   Using an estimate given by Voutier in 1996 (cf.
\cite{Vo}),
$$M(w)> \tau(d):=1+ \frac{1}{4}
\,\Big(\frac{\log\log d}{\log d}\Big)^3,
$$
which is the best known valid for every  $d>1$, not only
asymptotically, we get the following

\begin{teo}\label{entropy for infra}
Let $f: M \to M$ be a continuous map of a compact
infra-nilmanifold of dimension $d$. Then
\begin{itemize}
\item[a)] {either ${\bf h}(\phi_f)=0$,}
\item[b)] {or
${\bf h}(f)>  \log \tau(d)$.}
\end{itemize}
\end{teo}
\proof  Let $w(x)=w_1(x)\,\cdot w_2(x)\,\cdots\,w_k(x)$, $d_j=\deg
w_j$, $d_1\leq d_2\leq \cdots \leq d_k$, be a decomposition of the
characteristic polynomial of the linearization matrix $A_f$ into
irreducible terms. If
${\bf h}(\phi_f)>0$ then by Theorem B at least one eigenvalue of
$A_f$ has the absolute value larger than 1. Hence, by Theorem B,
using the property the sequence $\tau(n)$ is decreasing with respect to
$n$,
$$
{\bf h}(f)\geq \log {\underset{1\leq j\leq k}\prod} M(w_j(x)) \geq
\log M(w_{k}(x)) \geq \log \tau(d).
$$
\endproof

For other estimates of Mahler measure see for example \cite{EW}.


In particular from Smyth's theorem \cite{Sm} (which is a partial answer to
the Lehmer conjecture)
it follows

\begin{teo}\label{Smyth}
Let $f: M \to M$ be a continuous map of a compact
infra-nilmanifold of dimension $d$.
If the characteristic polynomial of $A_f$ is
non-reciprocal, i.e. the set of zeros is not invariant under
 the symmetry $\lambda\mapsto \lambda^{-1}$, and if ${\bf h}(\phi_f)>0$, then
$$
{\bf h}(f)\geq
{\underset{\lambda_j\in {\rm roots} \;w_j(x)}
\prod} \max(1,\vert \lambda_i^j\vert)
> \tau_0\,, $$ where $\tau_0$ is the real root of polynomial
$\tau^3-\tau-1$. \endproof
\end{teo}

One can check that the latter $\tau_0$ is greater
that $1.32471795$.
 In particular, $\tau_0$ does not depend neither on $w(x)$
nor on its degree $d$.

\

Note that ${\bf h}(f)\geq {\bf h}(\phi_f)$, $f\sim \phi_f$,
$A_f=A_{[f]}$ assert that Theorem \ref{Smyth} is a statement about
a homotopy property of $f$.  A special case is when $A_f$ is a
hyperbolic matrix invertible over integers, i.e. $\phi_f$ is an
Anosov automorphism, and $d$, the dimension of $M$, is odd. Then
obviously the characteristic polynomial of $A_f$ is
non-reciprocal, hence Theorem \ref{Smyth} applies and we obtain
$$
{\bf h}(f)\geq 1.32471795.
$$

This is in fact an easy case whose proof does not need the use of
Theorem B. Namely one can refer to Franks' theorem 
\cite[Theorem 2.2]{Franks}, 
saying that such a map $f$ is semiconjugate to
$\phi_f$, i.e. there exists a continuous map $\theta:M\to M$ such
that $\theta\circ f= \phi_f\circ\theta$. This $\theta$ is found to be homotopic to
identity, hence "onto". Therefore $ {\bf h}(f)\geq {\bf h}(\phi_f)$, see
Proposition \ref{entropy for covering}. It is easy to check that
if $f$ is an Anosov diffeomorphism then $A_f$ is a hyperbolic
invertible matrix.

\

{\bf Other remarks}

\

$\bullet$ The "projection -- shadowing" construction of $\Theta$ in the
proof of
Theorem B
in Section 3 and in \cite{MrPr} can be considered as a strengthening of Franks'
theorem to the case central direction exists.

\

$\bullet$ It is sufficient to assume $A_f$ is a hyperbolic endomorphism,  i.e. without
eigenvalues
of absolute value 1,  and without zero eigenvalues, to apply Franks'
theorem. Then $\phi_f$ is an Anosov endomorphism and the semiconjugacy
holds between the inverse limits, cf. \cite{Shub} and
\cite{Pr} .

\

$\bullet$ In the expanding case, i.e. if all the eigenvalues of $A_f$ have
absolute values larger than 1, the product is at least 2.
Therefore ${\bf h}(f)\geq \log 2$. In this case, instead of Theorem B, one
can refer to
Shub's theorem \cite{Shub} saying that $f$ is semiconjugate to $\phi_f$.

\

$\bullet$ Finally, if  $f$ itself is metric expanding on a compact
orientable
manifold (i.e. it expands all the distances between points close to each
other, at least by a constant factor
larger
than 1) or at least if $f$ is forward expansive, i.e. $\exists \delta>0$ such that
$\forall x\not=y$  $\exists n\ge 0$ with
$d(f^n(x),f^n(y))\ge \delta$ (as this implies expanding in an appropriate metric, see \cite[Section 3.6]{PU}),
 then for its degree $d(f)$ one has
immediately ${\bf h}(f)\ge \log |d(f)|\ge \log 2$, see \cite{Szl}.

Note that $f$ expanding (in a metric induced by a Riemannian metric) can
happen only on infra-nilmanifolds, \cite{Gr1}.

\

$\bullet$ In general ${\bf h}(f)\ge \log |d(f)|$  for all $f$ being $C^1$,
 see \cite{MiPr}.
 However the assumption that $f$ is $C^1$  is
essential in absence of the expanding
property, namely there are easy examples of continuous, but not
smooth maps $f$ for which ${\bf h}(f) < \log \vert \deg(f)\vert $.

\vskip 0.5cm \hskip 0.0cm  \vbox{\noindent
       Faculty of Mathematics and Comp. Sci.\newline
Adam Mickiewicz University of Pozna\'n\newline
 ul. Umultowska 87\newline
 61-614 Pozna{\'n}, Poland\newline
 e-mail: marzan@math.amu.edu.pl }

\vskip -2.5cm \hskip 8.5cm
 \vbox{\noindent Institute of Mathematics  
\newline
 Polish Academy of Sciences \newline
 ul. \'Sniadeckich 8 \newline
 00-950 Warszawa, Poland\newline
 e-mail: feliksp@impan.gov.pl
 }

\end{document}